\documentclass{amsart}
\usepackage{amssymb}
\usepackage{graphicx}
\usepackage[all]{xy}
\usepackage{a4wide}
\usepackage{hyperref}

\begin{document}
\newtheorem{theorem}{Theorem}
\newtheorem{corollary}[theorem]{Corollary}
\newtheorem{lemma}[theorem]{Lemma}
\theoremstyle{remark}
\newtheorem{remark}[theorem]{Remark}
\newcommand{\bbR}{\mathbb{R}}

\title{On The Remainder in the Taylor Theorem}%
\author{Lior Bary-Soroker}
\author{Eli Leher}
\address{Einstein Institute of Mathematics, The Hebrew University, Jerusalem 91904, Israel}%
\email{barylior@huji.math.ac.il}%
\address{School of Mathematical Sciences, Tel-Aviv University, Tel-Aviv 69978, Israel}%
\email{lehereli@post.tau.ac.il}%

\date{Dec. 24, 2008}
\begin{abstract}
We give a short straightforward proof for the bound of the reminder
term in the Taylor theorem. The proof uses only induction and the
fact that $f'\geq 0$ implies the monotonicity of $f$, so it might be
an attractive proof to give to undergraduate students.
\end{abstract}
\maketitle

\section{Introduction}
Let $f$ be an $n$-times differentiable function in a neighborhood of
$a\in \bbR$. Recall that the Taylor polynomial of order $n$ of $f$
at $a$ is the polynomial
\[
P_n (x) = f(a) + f'(a)(x-a) + \cdots + \frac{f^{(n)}(a)}{n!}(x-a)^n.
\]
It will be convenient to define $P_{-1}(x)=0$. Let $R_n = f - P_n$
be the remainder term. Then
\begin{theorem} [Lagrange's formula for the remainder]
If $f$ has an $(n+1)$th derivative in $[a,b]$ then there is some $a
\le \xi  \le b$ such that
\[
R_n(b)=\frac{f^{(n+1)}(\xi)}{(n+1)!}(b-a)^{n+1}.
\]
\end{theorem}

This formula is the main tool for bounding the remainder term of the
Taylor expansion in calculus classes, especially when this subject
is taught before integration. Therefore, one would like to have some
``natural" proof for it.
In \cite{Spivak} it is suggested that induction seems suitable, since $P_n'$ is the Taylor polynomial of
$f'$ of order $n-1$, hence $R_n'(x)$ is given by induction. The reason that this approach fails is that one cannot integrate $R_n'(x)$, since  $\xi=\xi(x)$ is implicit.

While we were teaching a first calculus course for chemistry and physics majors, we observed that this obstacle can be removed if we
slightly change the problem to finding a \textbf{bound} of the
remainder, which is all that is needed in order to show that the
Taylor series does in fact converge to the function.
From our personal experience, it seems that this approach enables
students to grasp the material more easily. Finally we mark that
Lagrange's formula can be deduced from the bound, as we show at the
end of this note.

The only fact needed in the proof is that a function with a positive
derivative is increasing. This can be easily proved with the mean value theorem
or without it (see \cite{Bers,Cohen}). As a direct
corollary one gets:

\begin{lemma}\label{lem:monotonicity}
Let $f,g$ be differentiable in a closed segment $[a,b]$. If
$f(a)=g(a)$ and $f'(x) \leq g'(x)$ for every $x \in (a,b)$, then
$f(x) \leq g(x)$ for every $a \leq x \le b$.
\end{lemma}

\begin{proof}
Take $h=g-f$. Thus $h' \geq 0$ and $h(a)=0$, hence $0\leq h(x)$,
i.e.\ $f(x)\leq g(x)$.
\end{proof}

\section{The main result}

\begin{theorem}  \label{Thm Bound}
Suppose that $f$ has an $(n+1)$th derivative in $[a,b]$ and that $m
\le f^{(n+1)}(x) \le M$ for every $a < x < b$. Then for any $a\leq
x\leq b$
\begin{equation}\label{eq:RemBnd}
\frac{m}{(n+1)!}(x-a)^{n+1} \le R_n(x) \le
\frac{M}{(n+1)!}(x-a)^{n+1}.
\end{equation}
\end{theorem}
\begin{proof}
By induction on $n$. For $n=-1$ the result is trivial.

For $n\geq 0$ write $f(x)=P_n(x)+R_n(x)$. Then
$f'(x)=P_n'(x)+R_n'(x)$. Note that $P_n'$ is the Taylor polynomial
of $f'$ of order $n-1$ and that $(f')^{(n)}=f^{(n+1)}$. Hence by
induction we have
\begin{equation}\label{eq:RemDer}
\frac{m}{n!}(x-a)^{n} \le R_n'(x) \le \frac{M}{n!}(x-a)^{n},
\end{equation}
for every $a \le x \le b$. Hence Lemma~\ref{lem:monotonicity} gives
the required inequality (since $\left(
\frac{(x-a)^{n+1}}{(n+1)!}\right)' = \frac{(x-a)^n}{n!}$).
\end{proof}

We conclude with a proof of Lagrange's classical formula. As
mentioned before, this might be omitted in calculus classes.

\begin{proof} Choose $m=\displaystyle\inf_{a\leq x\leq b}\{f^{(n+1)}(x)\}$ and
$M=\displaystyle\sup_{a\leq x\leq b}\{f^{(n+1)}(x)\}$ (if
$f^{(n+1)}$ is unbounded, we allow $m,M=\pm\infty$). Thus by
Theorem~\ref{Thm Bound} $R_n(b)=\frac{k}{(n+1)!}(b-a)^{n+1}$ for
some $m\leq k \leq M$. If one of equalities holds, then the result
is immediate from Theorem~\ref{Thm Bound}. Otherwise it follows
directly from Darboux's intermediate value theorem.
\end{proof}

\textbf{Acknowledgments:} We thank L.~Polterovich for his advices.


\begin{thebibliography}{99}
\bibitem{Bers}
L.~Bers, {\em On avoiding the mean value theorem}, Amer. Math.
Monthly {\bf 74} (1967) 583.
\bibitem{Cohen}
L.~W.~Cohen, {\em On being mean to the mean value theorem}, Amer.
Math. Monthly {\bf 74} 1967, 581--582.
\bibitem{Spivak}
M.~Spivak, {\em Calculus Third Edition}, Cambridge University Press,
(2006).
\end{thebibliography}
\end{document}